\documentclass[12 pt,twoside]{article}
\usepackage{graphicx}

\newtheorem{theorem}{\bf Theorem}[section]

\newtheorem{lemma}[theorem]{\bf Lemma}
\newtheorem{corollary}[theorem]{\bf Corollary}

\newtheorem{remark}[theorem]{\bf Remark}

\newtheorem{conjecture}{\bf Conjecture}

{{\footnotesize\sc }}
\input{amssym}

\begin{document}
\vspace{2 cm}
\title{{\Large Ramsey numbers for multiple copies of hypergraphs}\vspace{1 cm}}
\vspace{4cm}
\author{G.R. Omidi$^{\textrm{a},\textrm{c},1}$   \quad  G. Raeisi$^{\textrm{b},\textrm{c},2}$   \\[2pt]
{\small  $^{\textrm{a}}$Department of Mathematical Sciences, Isfahan University of Technology},\\
{\small Isfahan, 84156-83111, Iran}\\
{\small  $^{\textrm{b}}$Department of Mathematics, Shahrekord University,}\\
{\small  Shahrekord, P.O. Box 115, Iran}\\
{\small $^{\textrm{c}}$School of Mathematics, Institute for Research in Fundamental Sciences (IPM),}\\
{\small P.O. Box 19395-5746, Tehran, Iran }\\[2pt]
{\small romidi@cc.iut.ac.ir \quad g.raeisi@sci.sku.ac.ir, g.raeisi@math.iut.ac.ir}}

\date{}

\maketitle \footnotetext[1] {This research is partially
carried out in the IPM-Isfahan Branch and in part supported
by a grant from IPM (No. 91050416).} \vspace*{-0.5cm}

\maketitle \footnotetext[2] {This research was  in part supported
by a grant from IPM (No.91050018).} \vspace*{-0.5cm}

\footnotesize \medskip
\begin{abstract}\rm{}
\medskip

In this paper, for sufficiently large $n$ we determine the Ramsey
number $R(\mathcal{G}, n\mathcal{H})$ where $\mathcal{G}$ is a
$k$-uniform hypergraph  with the maximum independent set that
intersects each of the edges in $k-1$ vertices and $\mathcal{H}$
is a $k$-uniform hypergraph with a vertex so that the hypergraph
induced by the edges containing this vertex is a star. There are
several examples for such $\mathcal{G}$ and $\mathcal{H}$, among
them are any disjoint union of
 $k$-uniform hypergraphs involving loose paths, loose cycles,
tight paths, tight cycles with a multiple of $k$ edges, stars, Kneser hypergraphs and complete $k$-uniform $k$-partite hypergraphs for $\mathcal{G}$
 and linear hypergraphs for $\mathcal{H}$. As an application,
$R(m\mathcal{G},n\mathcal{H})$ is determined where $m$ or $n$ is
large and $\mathcal{G}$ and $\mathcal{H}$ are either loose paths, loose cycles,
tight paths, or stars. Also, $R(\mathcal{G},n\mathcal{H})$ is determined when  $\mathcal{G}$ is a bipartite graph with a matching saturating one of its color classes and $\mathcal{H}$ is an arbitrary graph for sufficiently large $n$.  Moreover, some bounds are given for
$R(m\mathcal{G},n\mathcal{H})$ which allow us to determine this
Ramsey number when $m\geq n$ and
$\mathcal{G}$ and $\mathcal{H}$ ($|V(\mathcal{G})|\geq |V(\mathcal{H})|$) are 3-uniform loose paths or cycles, $k$-uniform loose paths or cycles with at most 4 edges and $k$-uniform stars with 3 edges.
\vspace{.5cm}
 \\{ {Keywords}:{ \footnotesize Ramsey number, Hypergraph, Loose cycle, Loose path.\medskip}}
\end{abstract}
%%%%%%%%%%%%%%%%%%%%%%%%%%%%%%%%%%%%%%%%%%%%%%%%%%%%%%%%%%%%%%%%%%%%%%%%%%%%%%%%%%%%%%%%%%%%%%%%%%%%%
%%%%%%%%%%%%%%%%%%%%%%%%%%%%%%%%%%%%%%%%%%%%%%%%%%%%%%%%%%%%%%%%%%%%%%%%%%%%%%%%%%%%%%%%%%%%%%%%%%%%%%
%%%%%%%%%%%%%%%%%%%%%%%%%%%%%%%%%%%%%%%%%%%%%%%%%%%%%%%%%%%%%%%%%%%%%%%%%%%%%%%%%%%%%%%%%%%%%%%%%%%%%%
%%%%%%%%%%%%%%%%%%%%%%%%%%%%%%%%%%%%%%%%%%%%%%%%%%%%%%%%%%%%%%%%%%%%%%%%%%%%%%%%%%%%%%%%%%%%%%%%%%%%%%
%%%%%%%%%%%%%%%%%%%%%%%%%%%%%%%%%%%%%%%%%%%%%%%%%%%%%%%%%%%%%%%%%%%%%%%%%%%%%%%%%%%%%%%%%%%%%%%%%%%%%%
%%%%%%%%%%%%%%%%%%%%%%%%%%%%%%%%%%%%%%%%%%%%%%%%%%%%%%%%%%%%%%%%%%%%%%%%%%%%%%%%%%%%%%%%%%%%%%%%%%%%%%
\small
\newpage
\medskip
\section{\normalsize{Introduction}}

A {\it $k$-uniform hypergraph} $\mathcal{H}$ is a pair $(V,E)$
where $V$ is a set of vertices and $E$ is a set of $k$-subsets of
$V$ (the edges of $\mathcal{H}$). A hypergraph $\mathcal{H}$ is {\it linear} if the intersection
of every two edges of $\mathcal{H}$ has at most one element. As usual, the complete
$k$-uniform hypergraph on $p$ vertices is denoted by
$\mathcal{K}_p^k$ and for a given hypergraph $\mathcal{H}$,
$n\mathcal{H}$ is used to denote the $n$ disjoint copies of
$\mathcal{H}$. For a hypergraph $\mathcal{H}$, a set $S\subseteq
V(\mathcal{H})$ is called {\it independent} if there is no
edge of $\mathcal{H}$ contained in $S$. The {\it independence
number} of $\mathcal{H}$, denoted by $\alpha(\mathcal{H})$, is the
size of the greatest independent set in $\mathcal{H}$.

\bigskip
There are several natural definitions for a cycle and a path in
uniform hypergraphs. Those we focus on here are called {\it loose} and {\it tight}. By a
$k$-uniform {\it loose cycle} $\mathcal{C}_n^k$ (resp. {\it tight cycle} $\hat{\mathcal{C}}_n^k$), we mean the hypergraph with vertex set
$\{v_1,v_2,\ldots,v_{n(k-1)}\}$ (resp. $\{v_1,v_2,\ldots,v_{n}\}$) and with the set of $n$ edges
$e_i=\{v_1,v_2,\ldots, v_k\}+i(k-1)$, $i=0,1,\ldots, n-1$ (resp. $e_i=\{v_1,v_2,\ldots, v_k\}+i$, $i=0,1,\ldots, n-1$), where
we use mod $n(k-1)$ (resp. mod $n$) arithmetic and adding a number $t$ to a set
$S=\{v_1,v_2,\ldots, v_k\}$ means a shift, i.e., the set obtained
by adding $t$ to subscripts of each element of $S$. Similarly, a
$k$-uniform {\it loose path} $\mathcal{P}_n^k$ (resp. {\it tight path} $\hat{\mathcal{P}}_n^k$), is the hypergraph with vertex set
$\{v_1,v_2,\ldots,v_{n(k-1)+1}\}$ (resp. $\{v_1,v_2,\ldots,v_{n+k-1}\}$) and with the set of $n$ edges
$e_i=\{v_1,v_2,\ldots, v_k\}+i(k-1)$, $i=0,1,\ldots, n-1$ (resp. $e_i=\{v_1,v_2,\ldots, v_k\}+i$, $i=0,1,\ldots, n-1$). Also, by
a {\it star} $\mathcal{S}_n^k$ we mean the $k$-uniform hypergraph
with vertex set $\{v,v_1,v_2,\ldots,v_{n(k-1)}\}$ and with the set
of $n$ edges $e_i=\{v\}\cup(\{v_1,v_2,\ldots,
v_{k-1}\}+i(k-1))$, $i=0,1,\ldots, n-1$. For $k=2$ we get the
usall definitions of a cycle $C_n$, a path $P_n$ and a star
$K_{1,n}$ with $n$ edges.

\bigskip
For any given $k$-uniform hypergraphs $\mathcal{G}$ and
$\mathcal{H}$, the {\it Ramsey number} $R(\mathcal{G},
\mathcal{H})$ is the smallest positive integer $N$ such that in
every red-blue coloring of the edges of the complete $k$-uniform
hypergraph on $N$ vertices there is a monochromatic copy of
$\mathcal{G}$ in color red or a monochromatic copy of
$\mathcal{H}$ in color blue. The existence of such a positive
integer is guaranteed by Ramsey's classical result in
\cite{Ramsey}. The Ramsey number of graphs involving cycles and
paths are completely known (See \cite{FLPS,FS,Gerencser,Rosta}).
For the Ramsey number of tight paths and cycles the results in
\cite{HLPRRS} give the asymptotic behaviors of
$R(\hat{\mathcal{C}}^3_n, \hat{\mathcal{C}}^3_n)$ and
$R(\hat{\mathcal{P}}^3_n, \hat{\mathcal{P}}^3_n)$. The Ramsey
problem for loose paths and cycles were investigated by several
authors. It was proved in \cite{Ramsy number of loose cycle} that
$R(\mathcal{C}^3_n, \mathcal{C}^3_n)$ is asymptotically equal to
$\frac{5n}{2}$. Subsequently, Gy\'{a}rf\'{a}s {\it et. al.} in
\cite{Ramsy number of loose cycle for k-uniform} extended this
result to $k$-uniform loose cycles and proved that
$R(\mathcal{C}^k_n, \mathcal{C}^k_n)$ is asymptotically equal  to
$\frac{1}{2}(2k-1)n$. The proofs of all of these results are based
on the method of the Regularity Lemma. In \cite{prep}, the authors
determined the exact values of the Ramsey numbers of $k$-uniform
loose triangles and quadrangles.
\medskip

\begin{theorem}{\rm(\cite{prep})}\label{r(C3,C3)}
For every $k\geq 3$, $R(\mathcal{P}^k_3,\mathcal{P}^k_3)-1=R(\mathcal{C}^k_3,\mathcal{C}^k_3)=3k-2$ and
$R(\mathcal{P}^k_4,\mathcal{P}^k_4)-1=R(\mathcal{C}_4^k,\mathcal{C}_4^k)=4k-3$.
\end{theorem}

Also the Ramsey number of 3-uniform loose paths is determined when
one of the paths is significantly larger than the other. In the
other words, it is proved in \cite{M.O.R.S} that if $r\geq
\lfloor\frac{5s}{4}\rfloor$, then
$R(\mathcal{P}_r^3,\mathcal{P}_s^3)=2r+\lfloor\frac{s+1}{2}\rfloor.$
Recently in \cite{O.S}, the exact values of the Ramsey numbers of
$3$-uniform hypergraphs involving loose cycles and paths have been
determined as follows.
\medskip

\begin{theorem}{\rm(\cite{O.S})}\label{r(Cm,Cn)}
For every $r\geq s$,
$R(\mathcal{P}_r^3,\mathcal{P}_s^3)=R(\mathcal{C}_r^3,\mathcal{C}_s^3)+1=R(\mathcal{P}_r^3,\mathcal{C}_s^3)=2r+\Big\lfloor\frac{s+1}{2}\Big\rfloor$. Moreover,
$R(\mathcal{C}_r^3,\mathcal{P}_s^3)=2r+\Big\lfloor\frac{s-1}{2}\Big\rfloor$ if $r>s$.
\end{theorem}

The problem of determining the Ramsey numbers for multiple  copies
of graphs was first studied by Burr {\it et al.} in \cite{multiple
copies of graphs} where the authors found the Ramsey numbers for
multiple copies of  triangles and stars. More generally, Burr {\it
et al.} gave the following result on the Ramsey number of
connected graphs.

\medskip
\begin{theorem}{\rm(\cite{multiple copies of graphs})}\label{connectedgraphs}
Let $G$ and $H$ be connected graphs. Then there is a constant $c$,
depending only on $G$ and $H$, such that for sufficiently large
$n$
$$R(nG,nH)=(|V(G)|+|V(H)|-\min\{\alpha(G),\alpha(H)\})n+c.$$
\end{theorem}

In \cite{Burr 1}  and \cite{Burr 2}, Burr developed much more
powerful techniques to investigate the behavior of $R(nG,mH)$ when
either $m$  or $n$ is large. In particular, Burr proved  for fixed
$r,s\geq 4$ and sufficiently large $m$ or $n$, that

\begin{equation}\label{cycles}
R(mC_r,nC_s)=m|V(C_r)|+n|V(C_s)|-\min\{m\alpha(C_r),n\alpha(C_s)\}
- 1.
\end{equation}

\medskip
It is non-trivial to give a new version of Theorem \ref{connectedgraphs} for uniform hypergraphs,
here we do this by assuming some additional conditions. We give the exact value of $R(\mathcal{G},
n\mathcal{H})$ where $\mathcal{G}$ is a $k$-uniform hypergraph with the maximum independent set that intersects each of the edges in $k-1$ vertices, $\mathcal{H}$ has a vertex so that the hypergraph induced by the edges containing this vertex is a star and $n$ is sufficiently large. (See Theorem \ref{goodramsey}.)
Such evaluations are often possible in practice, as shown by
several examples; for instance, $R(m\mathcal{C}_r^k,n\mathcal{C}_s^k)$, $R(m\mathcal{C}_r^k,n\mathcal{P}_s^k)$,
$R(m\mathcal{P}_r^k,n\mathcal{P}_s^k)$, $R(m\hat{\mathcal{P}}_r^k,n\hat{\mathcal{P}}_s^k)$ and $R(m\mathcal{S}_r^k,n\mathcal{S}_s^k)$
are determined where $m$ or $n$ is
large and $k\geq 3$. (See Theorem \ref{nCycle}.) All these would satisfy a formula similar to (\ref{cycles}) with obvious substitutions. As an another example, $R(\mathcal{G},n\mathcal{H})$ is determined for a bipartite graph $\mathcal{G}$ with a matching saturating one of its color classes, an arbitrary graph $\mathcal{H}$ and sufficiently large $n$. Moreover, in Section 3, we give the exact values of Ramsey numbers for various
cases; for instance $R(m\mathcal{G},n\mathcal{H})$ is determined in the case when $m\geq n$ and
$\mathcal{G}$ and $\mathcal{H}$ ($|V(\mathcal{G})|\geq |V(\mathcal{H})|$) are 3-uniform loose paths or cycles, $k$-uniform loose paths or cycles with at most 4 edges and $k$-uniform stars with 3 edges.

\medskip
\section{\normalsize{$R(\mathcal{G},n\mathcal{H})$ for large $n$}}

 We  begin with some  definitions and notations.
Let $\mathcal{H}$ be a hypergraph and $S\subseteq V(\mathcal{H})$ (resp. $S\subseteq E(\mathcal{H})$).
 By {\it the induced hypergraph on $S$}, denoted by $\langle S\rangle$, we mean the hypergraph with
vertex set $S$ and those edges of $\mathcal{H}$ which are
contained in $S$ (resp. with vertex set $\bigcup_{e\in S} e$ and the edge set $S$). In the sequel, for a 2-edge coloring of
a uniform hypergraph $\mathcal{H}$, say red and blue, we denote by
$\mathcal{H}_{red}$ and $\mathcal{H}_{blue}$ the induced
hypergraph on edges of color red and blue, respectively.  A
{\it matching} in a hypergraph $\mathcal{H}$ is a set of mutually
disjoint edges and the {\it matching number}, $\nu(\mathcal{H})$, is
defined as  the size of the largest matching. A {\it covering} in
$\mathcal{H}$ is a set $S\subseteq V(\mathcal{H})$ such that any
edge of $\mathcal{H}$ intersects $S$. The {\it covering number} of
$\mathcal{H}$, $\tau(\mathcal{H})$, is defined as the size of the
smallest covering in $\mathcal{H}$.
%Also the {\it degree of a vertex} $v\in
%V(\mathcal{H})$ is  the number of edges of the largest star with center $v$.
%We use $\delta(\mathcal{H})$ to denote the {\it minimum degree} of
%vertices of $\mathcal{H}$.
%Let $\Omega$ be a set of $k$-uniform hypergraphs and $\mathcal{H}$ be an arbitrary hypergraph. We define $R(\Omega, \mathcal{H})$ to be the smallest positive %integer $N$ such that in every red-blue coloring of the
%edges of the complete $k$-uniform hypergraph $\mathcal{K}_N^k$
%there is a monochromatic copy of some member of $\Omega$ in color
%red or a monochromatic copy of $\mathcal{H}$ in color blue. Now
%define $\mathcal{D}(\mathcal{H})$ to be the set of all hypergraphs
%formed from the hypergraph $\mathcal{H}$ by removing a maximal set
%of independent vertices. For example,
%$\mathcal{D}(\mathcal{K}_n^r)=\{\mathcal{K}_{n-r+1}^r\}$.

\bigskip Another useful important variant in this paper is the
{\it strong independent set}. A {\it strong independent set} of a
$k$-uniform hypergraph $\mathcal{H}$ is an independent subset $S$
of vertices such that each edge of $\mathcal{H}$ has exactly one
vertex outside $S$. We denote by $\cal{F}$$_k$, the set of all
$k$-uniform hypergraphs which have a strong independent set. If
$\mathcal{H}\in \cal{F}$$_k$, the cardinality of the largest strong
independent set of $\mathcal{H}$ is called the {\it strong independence number} and is
denoted by $\alpha^\star(\mathcal{H})$. One can easily see that $\cal{F}$$_2$ is the set of all bipartite
graphs and for each $\mathcal{G}\in\cal{F}$$_2$, $\alpha^\star(\mathcal{G})$ is the size of the largest color class of $\mathcal{G}$ in all proper 2-colorings of $V(\mathcal{G})$. Clearly $\alpha^\star(\mathcal{H})\leq \alpha(\mathcal{H})$ for each $\mathcal{H}\in\cal{F}$$_k$.
A hypergraph $\mathcal{H}\in\cal{F}$$_k$ is called {\it good} if $\alpha^\star(\mathcal{H})=\alpha(\mathcal{H})$. We denote by $\cal{G}$$_k$, the set of all good $k$-uniform hypergraphs. Clearly $\cal{G}$$_k$ is closed under disjoint union, that is the disjoint union of every two hypergraphs in $\cal{G}$$_k$ is in $\cal{G}$$_k$. The following is a characterization of good uniform hypergraphs.

\begin{theorem}\label{goodhypergraphs}
Let $\mathcal{H}$ be a uniform hypergraph. Then $\mathcal{H}\in\cal{G}$$_k$ if and only if $V(\mathcal{H})$ can be partitioned into two subsets $V_1$ and $V_2$ so that each edge has one vertex in $V_1$ and $k-1$
vertices in $V_2$ and for each $S\subseteq V_1$,
$$\tau(\mathcal{H}_S)\geq |S|,$$
where $\mathcal{H}_S$ is the $k-1$-uniform hypergraph obtained from $\mathcal{H}$ by deleting
the vertices of $S$ from those edges of $\mathcal{H}$ which have nonempty intersection with $S$.
\end{theorem}

\noindent\textbf{Proof. } First, let $\mathcal{H}\in\cal{G}$$_k$. By the definition, the vertices of $\mathcal{H}$ can be partitioned into two subsets $V_1$ and $V_2$ so that each edge has one vertex in $V_1$ and $k-1$
vertices in $V_2$ and $|V_2|=\alpha^\star(\mathcal{H})=\alpha(\mathcal{H})$. Now, on contrary suppose that there is a set $S\subseteq V_1$ with
$|\tau(\mathcal{H}_S)|< |S|$ and assume that $S'\subseteq V(\mathcal{H}_S)$ is the minimum covering of $\mathcal{H}_S$. Clearly $V_2\cup S-S'$ is an independent set of $\mathcal{H}$ with more than $|V_2|=\alpha(\mathcal{H})$ vertices, a contradiction.

To prove the converse assume that $V(\mathcal{H})$ can be partitioned into two subsets $V_1$ and $V_2$ so that each edge has one vertex in $V_1$ and $k-1$
vertices in $V_2$ and for each $S\subseteq V_1$, we have
$\tau(\mathcal{H}_S)\geq |S|$. Suppose that $U\subseteq V(\mathcal{H})$ is the maximum independent set, $S=V_1\cap U$ and $S'= V(\mathcal{H}_S)-U$ is the covering of $\mathcal{H}_S$. Clearly $V_2\cup S-S'$ is an independent set of $\mathcal{H}$ and $U\subseteq V_2\cup S-S'$  and so $$|U|=|V_2\cup S-S'|=|V_2|+|S|-|S'|\leq|V_2|.$$ Since $V_2$ is an independent set we conclude that $|V_2|=\alpha(\mathcal{H})$ and since $|V_2|\leq \alpha^\star(\mathcal{H})\leq \alpha(\mathcal{H})$ we have $\mathcal{H}\in\cal{G}$$_k$.
$\hfill \blacksquare$

\bigskip A hypergraph is {\it $r$-regular} if each of its vertices lies on the $r$ edges. The following corollary shows that each regular element in  $\cal{F}$$_k$ is good.

\begin{corollary}\label{regularhypergraphs}
Let  $\mathcal{H}$ be $r$-regular and $\mathcal{H}\in\cal{F}$$_k$. Then $\mathcal{H}\in\cal{G}$$_k$.
\end{corollary}

\noindent\textbf{Proof. } Let $Y$ be the maximum strong independence set of $\mathcal{H}$ and  $X=V(\mathcal{H})-Y$. Hence each edge has one vertex in $X$ and $k-1$ vertices in $Y$. Now, let $S\subseteq X$. Clearly $\mathcal{H}_S$ has $r|S|$ edges and since $\mathcal{H}$ is $r$-regular the number of edges of $\mathcal{H}_S$ containing a vertex in $V(\mathcal{H}_S)$ is at most $r$. Hence $r\tau(\mathcal{H}_S)\geq r|S|$ and so the proof is complete by Theorem \ref {goodhypergraphs}.$\hfill \blacksquare$

\bigskip
For $k=2$ the condition $|\tau(\mathcal{H}_S)|\geq |S|$ in Theorem \ref{goodhypergraphs} is equivalent to the Hall condition for the existence a matching saturating $V_1$ in a bipartite graph $\mathcal{H}$. So we have the following result.

\bigskip
\begin{corollary}\label{bipartite-good}
$\mathcal{H}\in\cal{G}$$_2$ if and only if $\mathcal{H}$ is a bipartite graph with a matching saturating one of its color classes.
\end{corollary}

In the following we give some examples of good hypergraphs.

\medskip
\begin{remark}\end{remark}
Let $\mathcal{C}_n^k$ be the $k$-uniform loose cycle with
vertex set $\{v_1,v_2,\ldots,v_{n(k-1)}\}$ and the set of $n$
edges $e_i=\{v_1,v_2,\ldots, v_k\}+i(k-1)$ where $i=0,1,\ldots, n-1$.
It is clear to  see that if $n$ is even, then the set
$$S=\{v_1,v_{2k-1},v_{4k-3},\ldots,v_{(n-2)(k-1)+1}\},$$ is a
covering  and for odd $n$ the set
$$S=\{v_1,v_{2k-1},v_{4k-3},\ldots,v_{(n-3)(k-1)+1},v_{(n-1)(k-1)}\},$$ is a
covering for $\mathcal{C}_n^k$. In both cases, $T=V(\mathcal{C}_n^k)\setminus S$ with
$|T|=(k-1)n-\lfloor\frac{n+1}{2}\rfloor$ is a maximum strong
independent set (also a maximum independent set). Therefore
$\alpha^{\star}(\mathcal{C}_n^k)=\alpha(\mathcal{C}_n^k)=(k-1)n-\lfloor\frac{n+1}{2}\rfloor$.
By a similar argument, we have
$\alpha^{\star}(\mathcal{P}_n^k)=\alpha(\mathcal{P}_n^k)=(k-1)n-\lfloor\frac{n+1}{2}\rfloor+1$.
Also, one can easily see that $\alpha^{\star}(\hat{\mathcal{P}}_n^k)=\alpha(\hat{\mathcal{P}}_n^k)=n+k-2-\lfloor\frac{n-1}{k}\rfloor$ and for $n=kr$, $\alpha^{\star}(\hat{\mathcal{C}}_n^k)=\alpha(\hat{\mathcal{C}}_n^k)=n-r$. Thus, for every $k,n,r$
we have $\{\mathcal{C}^k_n,\mathcal{P}^k_n,\hat{\mathcal{C}}^k_{kr},\hat{\mathcal{P}}^k_n\}\subseteq \cal{G}$$_k$.

\bigskip
A $k$-uniform hypergraph {\it corresponding to a given graph $G$}, $\mathcal{H}_k(G)$, is a hypergraph
on $(k-2)|E(G)|+|V(G)|$ vertices and $|E(G)|$ edges where each of its edges is obtained by adding $k-2$ new vertices to an edge of $G$.
For example $\mathcal{P}_n^k=\mathcal{H}_k(P_n)$ and $\mathcal{S}_n^k=\mathcal{H}_k(K_{1,n})$. To see another family of good hypergraphs consider a tree $T$ with the property that the vertices of degree at least $3$ are independent. Since $\mathcal{P}^k_n\in\cal{G}$$_k$ one can easily see that $\mathcal{H}_k(T)\in\cal{G}$$_k$. Using Theorem \ref{goodhypergraphs}, we can show that for a bipartite graph $G$, $\mathcal{H}_k(G)\in\cal{G}$$_k$ if and only if $G\in\cal{G}$$_2$ or equivalently $G$ has a matching saturating one of its color classes. As an example, for a regular bipartite graph $G$ we have $\mathcal{H}_k(G)\in\cal{G}$$_k$.

\bigskip
A $k$-uniform hypergraph $\mathcal{H}$ is {\it $k$-partite} if its
vertices can be partitioned into $k$ classes such that each edge
intersects any class in exactly one vertex. A $k$-uniform
$k$-partite hypergraph is called {\it complete} if it contains all
possible edges. A complete $k$-uniform $k$-partite hypergraph with part
sizes $l_1,l_2,\ldots,l_k$ is denoted by
$\mathcal{K}(l_1,l_2,\ldots,l_k)$. For a complete $k$-uniform $k$-partite hypergraph, the vertices outside the smallest part is the maximum independent set and also the strong independent set and so for $l_1\geq l_2\geq\cdots\geq l_k$, $$\alpha^{\star}(\mathcal{K}(l_1,l_2,\ldots,l_k))=\alpha^{\star}(\mathcal{K}(l_1,l_2,\ldots,l_k))=\sum_{i=1}^{k-1} l_i.$$
Hence $\mathcal{K}(l_1,l_2,\ldots,l_k)\in \cal{G}$$_k$.

\bigskip
The {\it Kneser hypergraph} $\mathcal{KH}(n,r,k)$ is a $k$-uniform hypergraph whose vertices are the $r$-subsets of a given $n$-set $X$ and each of its edges contains $k$
mutually disjoint vertices. For $k=2$, this notion yields the usual definition of Kneser graphs. For $n=rk$ the edges of $\mathcal{KH}(n,r,k)$ containing those $r$-subsets of $X$
that do not contain a given element of $X$ is a strong independence set of $\mathcal{KH}(n,r,k)$ and so since $\mathcal{KH}(rk,r,k)$ is regular by Corollary \ref{regularhypergraphs}, $\mathcal{KH}(rk,r,k)$ is good. In this case $$\alpha(\mathcal{KH}(rk,r,k))=\alpha^{\star}(\mathcal{KH}(rk,r,k))={rk-1\choose r}.$$

\bigskip
\begin{theorem}\label{n large}
Let $\mathcal{G}\in \cal{F}$$_k$ and $\mathcal{H}$  be a $k$-uniform
hypergraph with a vertex so that the hypergraph induced by the edges containing this vertex is a star. Then for sufficiently large
$n$,
$$R(\mathcal{G},n\mathcal{H})\leq |V(\mathcal{G})|+n|V(\mathcal{H})|-\alpha^{\star}(\mathcal{G})-
1.$$
\end{theorem}
\noindent\textbf{Proof. } Assume that $\mathcal{H}$ has a vertex $v$ so that all edges containing $v$ makes a star with $\delta$ edges. Set $l=|V(\mathcal{H})|,
l^{\prime}=|V(\mathcal{G})|$,
$m=(k-1)(R(\mathcal{G},\mathcal{H})-1)l+\alpha^\star(\mathcal{G})$
and let $$n\geq
\Big(R(\mathcal{G},\mathcal{K}_m^k)-l^{\prime}+\alpha^{\star}(\mathcal{G})+1\Big)/l.$$
Assume that $p=l^{\prime}+nl-\alpha^{\star}(\mathcal{G})-1$ and consider a
2-edge colored $\mathcal{K}_p^k$ that contains no red copy of
$\mathcal{G}$. Set $V=V(\mathcal{K}_p^k)$. We will show that $\mathcal{K}_p^k$ must contain $n$ disjoint blue copies of
$\mathcal{H}$. First, we observe that $p\geq
R(\mathcal{G},\mathcal{K}^k_m)$ and so we have a blue copy of
$\mathcal{K}_m^k$ on  a set of vertices $B_1$. Find as many disjoint blue copies of $\mathcal{H}$ as possible in the induced hypergraph on $V\setminus B_1$,
 denoting the vertices  of these copies by $T_1$
and  $V\setminus (B_1\cup T_1)$ by $E_1$.  Clearly, $|E_1|\leq
R(\mathcal{G}, \mathcal{H})-1$, since the induced hypergraph on
$E_1$ does not contain a blue $\mathcal{H}$. If there is a vertex
$x\in E_1$ such that the degree of $x$ in $\langle B_1\cup \{x\}\rangle_{blue}$
is at least $\delta$, then this vertex and some $l-1$
vertices of $B_1$ span a blue $\mathcal{H}$. Transfer these $l$
vertices to $T_1$, and continue this process as long as possible.
This yields the three sets $E_2$, $B_2$ and $T_2$ such  that there is
no vertex $x\in E_2$ with degree at least $\delta$ in
$\langle B_2\cup\{x\}\rangle_{blue}$. Let $t=|E_1\setminus E_2|$. Clearly
$|E_2|\leq R(\mathcal{G}, \mathcal{H})-t-1$, $|B_2|= m-t(l-1)$ and
 the vertices of $T_2$ can be partitioned into the disjoint blue
copies of $H$. Now, a blue edge $e\in E(\mathcal{K}^k_p)$ is called {\it bad} if $|e\cap B_2|=k-1$ and $|e\cap E_2|=1$.
For a vertex $v\in E_2$, let $W_v$ be the set of all bad edges
containing $v$ and let $S_v$ be the largest star with center $v$ and
edges in $W_v$. Clearly by the condition on $E_2$, we have $|V(S_v)|\leq
(\delta-1)(k-1)+1$. Now transfer $V(S_v)\setminus
\{v\}$ into $T_2$ and for each $u\in E_2$, continue this process
as long as possible. This yields the sets $E_3$, $B_3$ and $T_3$.
Clearly, every  $[\frac{|T_3|}{l}]$$l$ vertices of $T_3$ can still be partitioned into disjoint blue
copies of $\mathcal{H}$ and for any $(k-1)$-set $S$ in $B_3$ and
for each vertex $v\in E_3$, the color of the $k$-set $S\cup \{v\}$
is red. On the other hand,
$$|B_3|\geq |B_2|-|E_2|(k-1)(\delta-1)\geq m-t(l-1)-(R(\mathcal{G},\mathcal{H})-t-1)(k-1)(\delta-1)\geq \alpha^\star (\mathcal{G}).$$
Therefore, $|E_3|<l^{\prime}-\alpha^\star (\mathcal{G})$ and so
$|B_3\cup T_3|\geq nl$. But then it is clear that the induced
hypergraph on $B_3\cup T_3$ contains $n$ disjoint blue copies of
$\mathcal{H}$. This observation completes the proof.

$\hfill \blacksquare$

\medskip Before giving some applications of Theorem \ref{n large},
we need the following lemma.

\medskip
\begin{lemma}\label{lower bound}
For every $k\geq 2$ and $k$-uniform hypergraphs $\mathcal{G}$ and
$\mathcal{H}$,
$$R(\mathcal{G},\mathcal{H})\geq \max\{|V(\mathcal{G})|+|V(\mathcal{H})|-\min\{\alpha(\mathcal{G}),\alpha(\mathcal{H})\}-1,
|V(\mathcal{G})|+\nu(\mathcal{H})-1,|V(\mathcal{H})|+\nu(\mathcal{G})-1\}.$$
\end{lemma}
\noindent\textbf{Proof. } First we exhibit
a 2-coloring, say red and blue, of the edges of the complete
$k$-uniform hypergraph on
$|V(\mathcal{G})|+|V(\mathcal{H})|-\min\{\alpha(\mathcal{G}),\alpha(\mathcal{H})\}-2$
vertices such that this coloring does not contain a red copy of
$\mathcal{G}$ and a blue copy of $\mathcal{H}$. For this purpose,
partition the vertex set into two parts $A$ and $B$, such that
$|A|=|V(\mathcal{G})|-\alpha(\mathcal{G})-1$ and
$|B|=|V(\mathcal{H})|-1$. We color all edges that contain a vertex
of $A$ red, and the rest blue. Now, this coloring can not contain
a blue copy of $\mathcal{H}$, since any such  copy must have all
vertices in $B$ and $|B|=|V(\mathcal{H})|-1$. Every red copy of
$\mathcal{G}$ would have to use $\alpha(\mathcal{G})+1$ vertices
of $B$, which is impossible since they would all be independent in
the red hypergraph. Thus $$R(\mathcal{G},\mathcal{H})\geq
|V(\mathcal{G})|+|V(\mathcal{H})|-\alpha(\mathcal{G})-1.$$ By
symmetry $$R(\mathcal{G},\mathcal{H})\geq
|V(\mathcal{G})|+|V(\mathcal{H})|-\alpha(\mathcal{H})-1.$$
Combining the two inequalities yields the desired result.

\medskip
Now partition the vertex set of the
complete $k$-uniform hypergraph on
$|V(\mathcal{G})|+\nu(\mathcal{H})-2$ vertices into two parts $A$
and $B$, such that $|A|=\nu(\mathcal{H})-1$ and
$|B|=|V(\mathcal{G})|-1$. We color all edges that contain a vertex
of $A$ blue, and the rest red. Now, this coloring can not contain
a red copy of $\mathcal{G}$, since any such  copy must have all
vertices in $B$ and $|B|=|V(\mathcal{G})|-1$. Also the maching
number of the blue hypergraph is at most $\nu(\mathcal{H})-1$ and
so it can not contain a copy of $\mathcal{H}$. Hence $R(\mathcal{G},\mathcal{H})\geq |V(\mathcal{G})|+\nu(\mathcal{H})-1$ and by
symmetry again $R(\mathcal{G},\mathcal{H})\geq|V(\mathcal{H})|+\nu(\mathcal{G})-1$.

$\hfill \blacksquare$
%%%%%%%%%%%%%%%%%%%%%%%%%%%%%%%%%%%%%%%%%%%%%%%%%%%%%%%%%%%%%%%%%%%%%%%%%%%%%%%%%%%%%%%%%%%%%%%%%%%%%%
%%%%%%%%%%%%%%%%%%%%%%%%%%%%%%%%%%%%%%%%%%%%%%%%%%%%%%%%%%%%%%%%%%%%%%%%%%%%%%%%%%%%%%%%%%%%%%%%%%%%%%

The following result is an immediate consequence of Theorem
\ref{n large} and Lemma \ref{lower bound}.

\bigskip
\begin{theorem}\label{goodramsey}
Assume that $\mathcal{G}\in \cal{G}$$_k$ and $\mathcal{H}$ is a $k$-uniform hypergraph with a vertex so that the hypergraph induced by the edges containing this vertex is a star. Then for sufficiently large $n$,
$$R(\mathcal{G},n\mathcal{H})=|V(\mathcal{G})|+n|V(\mathcal{H})|-\alpha(\mathcal{G})-
1.$$
\end{theorem}

\bigskip
\begin{corollary}
Assume that $\mathcal{G}\in \cal{G}$$_k$ and $\mathcal{H}$ is a $k$-uniform linear hypergraph. Then for sufficiently large $n$,
$$R(\mathcal{G},n\mathcal{H})=|V(\mathcal{G})|+n|V(\mathcal{H})|-\alpha(\mathcal{G})-
1.$$
\end{corollary}

Note that for $\mathcal{G}\in \cal{F}$$_k$, since  $\alpha^\star(\mathcal{G})\leq \alpha(\mathcal{G})\leq |V(\mathcal{G})|-\nu(\mathcal{G})$, the condition  $\alpha^\star(\mathcal{G})= |V(\mathcal{G})|-\nu(\mathcal{G})$ implies $\mathcal{G}\in \cal{G}$$_k$ and so we do not have any new result if we consider this condition instead of $\mathcal{G}\in \cal{G}$$_k$ in Theorem \ref{goodramsey}. Using Remark 2.4 and Theorem \ref{goodramsey}, we have the following results.

%$$R(m\hat{\mathcal{C}}_{kr}^k,n\hat{\mathcal{C}}_{ks}^k)=krm+ksn-\min\{m\alpha(\hat{\mathcal{C}}_{kr}^k),n\alpha(\hat{\mathcal{C}}_{ks}^k)\}-1,$$
%$$R(m\mathcal{KH}(rk,r,k),n\mathcal{KH}(sk,s,k))=m{rk\choose r}+n{sk\choose s}-\min\{m{rk-1\choose r},n{sk-1\choose s}\}-1.$$
%\\{\rm(ii)}~ Let $1\leq l_1\leq l_2\leq \cdots \leq l_k$, $1\leq t_1\leq t_2\leq \cdots \leq t_k$, $r=\sum_{i}l_i$ and $s=\sum_{i}t_i$. Then
%$$R(m\mathcal{K}(l_1,
%l_2,\ldots, l_k),n\mathcal{K}(t_1, t_2,\ldots, t_k))=mr+ns-\min\{m(r-l_1),n(s-t_1)\}-1.$$

\bigskip
\begin{theorem}\label{nCycle}
Assume that either $m$ or $n$ is sufficiently large, $k\geq 3$ and $\mathcal{G}$, $\mathcal{H}\in \cal{G}$$_2$. Then
$$R(m\mathcal{C}_r^k,n\mathcal{C}_s^k)=(k-1)rm+(k-1)sn-\min\{m\alpha(\mathcal{C}_r^k),n\alpha(\mathcal{C}_s^k)\}-1,$$
$$R(m\mathcal{P}_r^k,n\mathcal{C}_s^k)=((k-1)r+1)m+(k-1)sn-\min\{m\alpha(\mathcal{P}_r^k),n\alpha(\mathcal{C}_s^k)\}-1,$$
$$R(m\mathcal{P}_r^k,n\mathcal{P}_s^k)=((k-1)r+1)m+((k-1)s+1)n-\min\{m\alpha(\mathcal{P}_r^k),n\alpha(\mathcal{P}_s^k)\}-1,$$
$$R(m\hat{\mathcal{P}}_r^k,n\hat{\mathcal{P}}_s^k)=(r+k-1)m+(s+k-1)n-\min\{m\alpha(\hat{\mathcal{P}}_r^k),n\alpha(\hat{\mathcal{P}}_s^k)\}-1,$$
$$R(m\mathcal{S}_r^k,n\mathcal{S}_s^k)=(k-1)(mr+ns)+m+n-\min\{mr(k-1),ns(k-1)\}-1,$$
$$R(m\mathcal{G},n\mathcal{H})= m|V(\mathcal{G})|+n|V(\mathcal{H})|-\min\{m\alpha(\mathcal{G}),n\alpha(\mathcal{H})\}-1.$$
\end{theorem}

\begin{theorem}\label{bipartite}
Assume that $\mathcal{G}$ is a bipartite graph with a matching saturating one of its color classes and $\mathcal{H}$ is an arbitrary graph. Then for sufficiently large $n$,
$$R(\mathcal{G},n\mathcal{H})=n|V(\mathcal{H})|+\nu(\mathcal{G})-
1.$$
\end{theorem}

\bigskip
\section{\normalsize  Multiple copies of loose paths and cycles}

In this section, we provide the exact values of
$R(m\mathcal{G},n\mathcal{H})$ for every $m\geq n$ and particular
hypergraphs  $\mathcal{G}$ and $\mathcal{H}$ with $\alpha(\mathcal{G})\geq \alpha(\mathcal{H})$, for example,
$k$-uniform loose triangles, loose quadrangles, stars with maximum degree 3  and 3-uniform loose paths and cycles. Before that, we need
the following.

\medskip
\begin{lemma}\label{lower bound1}
For every $k\geq 2$ and $k$-uniform hypergraphs $\mathcal{G}$,
$\mathcal{H}$ and $\mathcal{F}$,
$$R(\mathcal{G},\mathcal{H}\cup \mathcal{F})\leq \max\{R(\mathcal{G},\mathcal{F})+|V(\mathcal{H})|,R(\mathcal{G},\mathcal{H})\},$$
$$R(m\mathcal{G},n\mathcal{H})\leq R(\mathcal{G},\mathcal{H})+(m-1)|V(\mathcal{G})|+(n-1)|V(\mathcal{H})|.$$
\end{lemma}
\noindent\textbf{Proof. }Let
$r=\max\{R(\mathcal{G},\mathcal{F})+|V(\mathcal{H})|,R(\mathcal{G},\mathcal{H})\}$
and the edges of $K=\mathcal{K}^k_{r}$ be 2-colored red and blue.
If there is no red $\mathcal{G}$, then there must certainly be a
blue $\mathcal{H}$. Remove the vertices of this blue copy of
$\mathcal{H}$ from $K$. Among the remaining vertices there
must be a blue $\mathcal{F}$ or a red $\mathcal{G}$. Hence $K$
contains either a red $\mathcal{G}$ or a blue $\mathcal{H}\cup
\mathcal{F}$, and the first inequality follows. The second
inequality follows by applying the first
inequality. $\hfill \blacksquare$
%%%%%%%%%%%%%%%%%%%%%%%%%%%%%%%%%%%%%%%%%%%%%%%%%%%%%%%%%%%%%%%%%%%%%%%%%%%%%%%%%%%%%%%%%%%%%%%%%%%%%%
%%%%%%%%%%%%%%%%%%%%%%%%%%%%%%%%%%%%%%%%%%%%%%%%%%%%%%%%%%%%%%%%%%%%%%%%%%%%%%%%%%%%%%%%%%%%%%%%%%%%%%

\medskip
\begin{theorem}\label{main theorem1}
Assume that $k\geq 2$, $m\geq n\geq 2$ and  $\mathcal{G}$ and
$\mathcal{H}$ are arbitrary $k$-uniform hypergraphs. Then
$$R(n\mathcal{G},m\mathcal{H})\leq R((n-1)\mathcal{G},(m-1)\mathcal{H})+R(\mathcal{G},\mathcal{H})+1.$$
In particular,
$$R(n\mathcal{G},m\mathcal{H})\leq R(\mathcal{G},(m-n+1)\mathcal{H})+(n-1)(R(\mathcal{G},\mathcal{H})+1).$$
\end{theorem}
\noindent\textbf{Proof. }Set
$t=R((n-1)\mathcal{G},(m-1)\mathcal{H})+R(\mathcal{G},\mathcal{H})+1$
and let $K=\mathcal{K}_{t}^k$ be 2-edge colored red and blue. We
find either a red $n\mathcal{G}$ or a blue $m\mathcal{H}$. We have $t\geq R(\mathcal{G},\mathcal{H})$ and thus we
may assume that $K$ contains a red copy of $\mathcal{G}$ (we have the same proof if $K$ contains a blue copy of $\mathcal{H}$). Discard
this copy. Since the number of remaining vertices is greater than
$R((n-1)\mathcal{G},(m-1)\mathcal{H})$, there is either a red copy
of $(n-1)\mathcal{G}$ or a blue copy of $(m-1)\mathcal{H}$. In the
first case, we have a red $n\mathcal{G}$ and so  we are done. Thus
we may assume that $K$ contains a blue $(m-1)\mathcal{H}$ and
therefore we have a red copy of $\mathcal{G}$ and a blue copy of
$\mathcal{H}$. Among red-blue copies of $\mathcal{G}$ and
$\mathcal{H}$ choose red-blue copies with maximum intersection.
Let $\mathcal{G}^{\prime}$ and $\mathcal{H}^{\prime}$ be such
copies. We must have $p=|V(\mathcal{G}^{\prime}\cup
\mathcal{H}^{\prime})|\leq R(\mathcal{G},\mathcal{H})+1$. Indeed,
let $p\geq R(\mathcal{G},\mathcal{H})+2$. Clearly $p\geq
\max\{|V(\mathcal{G}^{\prime})|,|V(\mathcal{H}^{\prime})|\}+2$ and
hence $V(\mathcal{G}^{\prime})\setminus V(\mathcal{H}^{\prime})$
and $V(\mathcal{H}^{\prime})\setminus V(\mathcal{G}^{\prime})$ are
non-empty. Choose $v_1\in V(\mathcal{G}^{\prime})\setminus
V(\mathcal{H}^{\prime})$ and $v_2\in
V(\mathcal{H}^{\prime})\setminus V(\mathcal{G}^{\prime})$. Set
$U=(V(\mathcal{G}^{\prime})\cup
V(\mathcal{H}^{\prime}))\setminus \{v_1,v_2\}$. Since $|U|=p-2\geq R(\mathcal{G},\mathcal{H})$, we have either a
red $\mathcal{G}$ or a blue $\mathcal{H}$, say $\mathcal{F}$. If
$\mathcal{F}$ is red, then $|\mathcal{F}\cap
\mathcal{H}^{\prime}|> |\mathcal{G}^{\prime}\cap
\mathcal{H}^{\prime}|$. If $\mathcal{F}$ is blue, then
$|\mathcal{F}\cap \mathcal{G}^{\prime}|> |\mathcal{G}^{\prime}\cap
\mathcal{H}^{\prime}|$. Both cases, contradict the choice of
$\mathcal{G}^{\prime}$ and $\mathcal{H}^{\prime}$. Therefore
$|V(\mathcal{G}^{\prime}\cap \mathcal{H}^{\prime})|\leq
R(\mathcal{G},\mathcal{H})+1$. Remove the vertices of $\mathcal{G}^{\prime}\cup
\mathcal{H}^{\prime}$ from $K$. The hypergraph on the remaining
vertices contains either a red $(n-1)\mathcal{G}$ or a blue
$(m-1)\mathcal{H}$, say $T$, to which we add the appropriately
colored copy of $\mathcal{G}$ and $\mathcal{H}$ to obtain a red
$n\mathcal{G}$ or a blue $m\mathcal{H}$. This observation completes the
proof of the first inequality. The second inequality follows from
repeated application of the first inequality, which completes the
proof.

$\hfill \blacksquare$
%%%%%%%%%%%%%%%%%%%%%%%%%%%%%%%%%%%%%%%%%%%%%%%%%%%%%%%%%%%%%%%%%%%%%%%%%%%%%%%%%%%%%%%%%%%%%%%%%%%%%
%%%%%%%%%%%%%%%%%%%%%%%%%%%%%%%%%%%%%%%%%%%%%%%%%%%%%%%%%%%%%%%%%%%%%%%%%%%%%%%%%%%%%%%%%%%%%%%%%%%%%%

\bigskip
As an easy, but useful application of Theorem \ref{main theorem1},
we have the following corollary.

\medskip
\begin{corollary}\label{Ramsey nH,nG}
Assume that $k\geq 2$, $n\geq 1$ and  $\mathcal{G}$ and
$\mathcal{H}$ are $k$-uniform hypergraphs. Then
$$R(n\mathcal{G},n\mathcal{H})\leq (R(\mathcal{G},\mathcal{H})+1)n-1.$$
\end{corollary}
%\noindent\textbf{Proof. }

%%%%%%%%%%%%%%%%%%%%%%%%%%%%%%%%%%%%%%%%%%%%%%%%%%%%%%%%%%%%%%%%%%%%%%%%%%%%%%%%%%%%%%%%%%%%%%%%%%%%%
%%%%%%%%%%%%%%%%%%%%%%%%%%%%%%%%%%%%%%%%%%%%%%%%%%%%%%%%%%%%%%%%%%%%%%%%%%%%%%%%%%%%%%%%%%%%%%%%%%%%%%

\medskip
Using Lemma \ref{lower bound} and Corollary \ref{Ramsey nH,nG}, we
get the following theorem.

\medskip
\begin{theorem}\label{main theorem}
Assume that $k\geq 2$ and $m\geq n\geq 1$ are positive integers
and $\mathcal{G}$ and $\mathcal{H}$ are $k$-uniform hypergraphs with
$\alpha(\mathcal{G})\geq \alpha(\mathcal{H})$ and
$R(\mathcal{G},\mathcal{H})=
|V(\mathcal{G})|+|V(\mathcal{H})|-\alpha(\mathcal{H})-1$. Then

$$R(m\mathcal{G},n\mathcal{H})=m|V(\mathcal{G})|+n|V(\mathcal{H})|-n\alpha(\mathcal{H})-1.$$

\end{theorem}
\noindent\textbf{Proof. }Let
$t=m|V(\mathcal{G})|+n|V(\mathcal{H})|-n\alpha(\mathcal{H})-1$ and
$K=\mathcal{K}_{t}^k$ be 2-edge colored red and blue. We use
induction on $m$ to prove that either $m\mathcal{G}\subseteq
K_{red}$ or $n\mathcal{H}\subseteq K_{blue}$.
Clearly for $m=1$, the result is true and so we may assume that
$m\geq 2$. For $m=n$, the result follows from Corollary
\ref{Ramsey nH,nG} and so we may assume that $m-1\geq n$. Since by the
induction hypothesis $$R((m-1)\mathcal{G},n\mathcal{H})\leq
(m-1)|V(\mathcal{G})|+n|V(\mathcal{H})|-n\alpha(\mathcal{H})-1<t,$$
we may assume that we have $(m-1)\mathcal{G}\subseteq
K_{red}$, otherwise we can find $n$
disjoint blue copies of $\mathcal{H}$. Now, remove the vertices of
a red $\mathcal{G}$ from $K$ and use the induction hypothesis to the
coloring on the remaining
$(m-1)|V(\mathcal{G})|+n|V(\mathcal{H})|-n\alpha(\mathcal{H})-1$
vertices to find either a $(m-1)\mathcal{G}\subseteq
K_{red}$ or a $n\mathcal{H}\subseteq
K_{blue}$. If $n\mathcal{H}\subseteq
K_{blue}$ we are done, otherwise $(m-1)\mathcal{G}\subseteq
K_{red}$, adding the deleted red colored $\mathcal{G}$
to the red $(m-1)\mathcal{G}$, we obtain $m$
disjoint red copies of $\mathcal{G}$, which shows that
$R(m\mathcal{G},n\mathcal{H})\leq t$.
\medskip

For the lower bound, employ Lemma \ref{lower bound}, where
$\mathcal{G}$ is replaced by $m\mathcal{G}$ and $\mathcal{H}$
is replaced by $n\mathcal{H}$, which completes the proof.
 $\hfill
\blacksquare$
%%%%%%%%%%%%%%%%%%%%%%%%%%%%%%%%%%%%%%%%%%%%%%%%%%%%%%%%%%%%%%%%%%%%%%%%%%%%%%%%%%%%%%%%%%%%%%%%%%%%%%
%%%%%%%%%%%%%%%%%%%%%%%%%%%%%%%%%%%%%%%%%%%%%%%%%%%%%%%%%%%%%%%%%%%%%%%%%%%%%%%%%%%%%%%%%%%%%%%%%%%%

\bigskip
In the rest of this section, we use Theorem \ref{main theorem} to
give some corollaries. Before that, we complete the determining of the Ramsey numbers of loose paths and cycles with at most 4 edges.

\medskip
\begin{lemma}\label{star-cycle}
For every $k\geq 3$, $R(\mathcal{C}_3^k,\mathcal{C}_4^k)=4k-3$ and
$R(\mathcal{S}_3^k,\mathcal{S}_3^k)=3k-2$.
\end{lemma}
\noindent\textbf{Proof. }Using Lemma \ref{lower bound},
$R(\mathcal{C}_3^k,\mathcal{C}_4^k)\geq4k-3$ and
$R(\mathcal{S}_3^k,\mathcal{S}_3^k)\geq 3k-2$.  To prove that
$R(\mathcal{C}^k_3,\mathcal{C}^k_4)\leq 4k-3$, suppose that the
edges of $K=\mathcal{K}_{4k-3}^k$ are arbitrary colored red and
blue. We prove that $K$ contains a red copy of $\mathcal{C}_3^k$
or a blue copy of $\mathcal{C}_4^k$. Since
$R(\mathcal{C}_4^k,\mathcal{C}_4^k)=4k-3$, we may assume that $K$
contains a red copy of $\mathcal{C}_4^k$. Let
$e_i=\{v_1,v_2,\ldots, v_k\}+i(k-1)$ mod $4$, $0\leq i\leq3$, be
the edges of $\mathcal{C}_4^k\subseteq K_{red}$ and $v$ be the
remaining vertex which is not covered by this copy of
$\mathcal{C}_4^k$. Set $e_{0}^{\prime}=\{v_1,v_{2k-1},\ldots,
v_{3k-3}\}$, $e_{1}^{\prime}=\{v_{3k-3},v_{4k-4},v_2,\ldots,
v_{k-2},v\}$,
$e_{2}^{\prime}=\{v_{2k-2},v_{3k-2},v_{3k},v_{3k+1},\ldots,
v_{4k-4},v_{k}\}$ and
$e_{3}^{\prime}=\{v_{3k-1},v_{k-1},v_{k+2}\ldots, v_{2k-1}\}$. If
one of $e_i^{\prime}$ is red, we have a red copy of
$\mathcal{C}_3^k$, otherwise
$e_0^{\prime},e_1^{\prime},e_2^{\prime},e_3^{\prime}$ form a blue
copy of $\mathcal{C}_4^k$ which shows that
$R(\mathcal{C}_3^k,\mathcal{C}_4^k)\leq 4k-3$. To see
$R(\mathcal{S}_3^k,\mathcal{S}_3^k)\leq 3k-2$, let  the edges of
$K=\mathcal{K}_{3k-2}^k$ be arbitrary colored red and blue. By
Theorem \ref{r(C3,C3)}, we have a monochromatic, say red, copy of
$\mathcal{C}^k_3$. Assume that $e_1=\{v,v_1,\ldots, v_{k-2},u\}$,
$e_2=\{u,u_1,\ldots, u_{k-2},w\}$ and $e_3=\{w,w_1,\ldots
,w_{k-2},v\}$ are the edges of this copy of $\mathcal{C}_3^k$ and
$T=\{t\}$ be the remaining vertex of $K$. If one of the edges
$e_1^{\prime}=\{t,w,v_1,\ldots, v_{k-2}\}$,
$e_2^{\prime}=\{t,u,w_1,\ldots, w_{k-2}\}$ or
$e_3^{\prime}=\{t,v,u_1,\ldots, u_{k-2}\}$ is red, then we have a
red copy of $\mathcal{S}_3^k\subseteq K_{red}$, otherwise
$e_1^{\prime}e_2^{\prime}e_3^{\prime}$ form a
$\mathcal{S}_3^k\subseteq K_{blue}$. This observation completes
the proof.
 $\hfill\blacksquare$
%%%%%%%%%%%%%%%%%%%%%%%%%%%%%%%%%%%%%%%%%%%%%%%%%%%%%%%%%%%%%%%%%%%%%%%%%%%%%%%%%%%%%%%%%%%%%%%%%%%%%%
%%%%%%%%%%%%%%%%%%%%%%%%%%%%%%%%%%%%%%%%%%%%%%%%%%%%%%%%%%%%%%%%%%%%%%%%%%%%%%%%%%%%%%%%%%%%%%%%%%%%

%
%\medskip
\medskip
\begin{lemma}{\rm(\cite{prep})}\label{lower bound}
Let $n\geq m\geq 3$, $k\geq 3$ and $r=(k-1)n+\lfloor\frac{m+1}{2}\rfloor$. Then
\\{\rm(i)} $R(\mathcal{C}^k_n,\mathcal{C}^k_m)\geq
r-1$ and also $r$ is a lower bound for both
$R(\mathcal{P}^k_n,\mathcal{P}^k_m)$ and
$R(\mathcal{P}^k_n,\mathcal{C}^k_m)$.
\\{\rm(ii)} Assume that $\mathcal{K}^k_{r}$ is 2-edge
colored red and blue. If $\mathcal{C}_{n}^k\subseteq
\mathcal{F}_{red}$, then either $\mathcal{P}_n^k\subseteq
\mathcal{F}_{red}$ or $\mathcal{P}_m^k\subseteq
\mathcal{F}_{blue}$. Also, if $\mathcal{C}_{n}^k\subseteq
\mathcal{F}_{red}$, then either $\mathcal{P}_n^k\subseteq
\mathcal{F}_{red}$ or $\mathcal{C}_m^k\subseteq
\mathcal{F}_{blue}$.
\end{lemma}

The following theorem is a direct consequence of Lemmas \ref{star-cycle} and \ref{lower bound} and Theorem \ref{r(C3,C3)}.

\begin{lemma}\label{smallcase}
For every $k\geq 3$, we have
$$R(\mathcal{P}^k_3,\mathcal{P}^k_4)=R(\mathcal{C}^k_3,\mathcal{P}^k_4)=4k-2,$$
$$R(\mathcal{P}^k_3,\mathcal{C}^k_4)=4k-3,~~~~~R(\mathcal{C}^k_3,\mathcal{P}^k_3)=3k-1.$$
\end{lemma}

\bigskip Now, using the known result $R(P_r,P_s)=r+\lfloor\frac{s}{2}\rfloor-1$ for $r\geq s$ due to Gerencs\'{e}r and Gy\'{a}rf\'{a}s in \cite{Gerencser} and Theorems \ref{r(C3,C3)}, \ref{r(Cm,Cn)} and \ref{main theorem} and Lemmas
\ref{star-cycle} and \ref{smallcase} , we have the following theorems.
\bigskip
\begin{theorem}\label{ramsey number for k-psths}
If $m\geq n\geq 1$ and $k\geq 3$, then
\\{\rm(i)}~$R(m\mathcal{C}_3^k,n\mathcal{C}_3^k)=m(3k-3)+2n-1$ and $R(m\mathcal{C}_4^k,n\mathcal{C}_3^k)=R(m\mathcal{C}_4^k,n\mathcal{C}_4^k)=m(4k-4)+2n-1,$
\\{\rm(ii)}~$R(m\mathcal{P}_3^k,n\mathcal{C}_3^k)=m(3k-2)+2n-1$ and $R(m\mathcal{P}_4^k,n\mathcal{C}_3^k)=R(m\mathcal{P}_4^k,n\mathcal{C}_4^k)=m(4k-3)+2n-1,$
\\{\rm(iii)}~$R(m\mathcal{P}_3^k,n\mathcal{P}_3^k)=m(3k-2)+2n-1$ and $R(m\mathcal{P}_4^k,n\mathcal{P}_3^k)=R(m\mathcal{P}_4^k,n\mathcal{P}_4^k)=m(4k-3)+2n-1,$
\\{\rm(iv)}~$R(m\mathcal{S}^k_3,n\mathcal{S}^k_3)=m(3k-2)+n-1,$
\\{\rm(v)}~For every $k$-uniform hypergraph $\mathcal{H}$,
$R(m\mathcal{H},n\mathcal{K}_k^k)=m|V(\mathcal{H})|+n-1.$ In
particular $R(m\mathcal{K}_k^k,n\mathcal{K}_k^k)=mk+n-1$.
\end{theorem}

\begin{theorem}\label{ramsey number for multiple 3-paths}
For every $m\geq n\geq 1$ and $r\geq s\geq 1$ we have the
following.\\{\rm(i)}~$R(m\mathcal{P}_r^3,n\mathcal{P}_s^3)=R(m\mathcal{P}_r^3,n\mathcal{C}_s^3)=(2r+1)m+\lfloor\frac{s+1}{2}\rfloor n-1,$
\\{\rm(ii)}~$R(m\mathcal{C}_r^3,n\mathcal{C}_s^3)=2rm+\lfloor\frac{s+1}{2}\rfloor n-1,$
\\{\rm(iii)}~$R(m\mathcal{C}_r^3,n\mathcal{P}_s^3)=2rm+\lfloor\frac{s+1}{2}\rfloor n-1$ if $r>s$,
\\{\rm(iv)}~$R(mP_r,nP_s)=rm+
\lfloor\frac{s}{2}\rfloor n-1$.
\end{theorem}

\section{\normalsize Concluding remarks}

By Lemma \ref{lower bound}, for every $k\geq 3$,
$m\geq n$ and $r\geq s\geq 3$,
\begin{eqnarray}
R(m\mathcal{C}_r^k,n\mathcal{C}_s^k)\geq
(k-1)rm+\Big\lfloor\frac{s+1}{2}\Big\rfloor n-1.
\end{eqnarray}
By Theorem \ref{nCycle}, we have equality in $(2)$ if $m$ is
sufficiently large. It would be interesting to decide whether this
natural lower bound is always the exact value of the Ramsey number. The case
$k=3$ follows from Theorem \ref{ramsey number for multiple
3-paths}. Based on these observations and with the same discussions for the Ramsey numbers of multiple copies of hypergraphs involving loose paths, loose cycles,
tight paths and tight cycles we pose the following
conjecture.

\medskip
\begin{conjecture}
For every $k\geq 3$, $m\geq n$ and $r\geq s\geq 3$,
$$R(m\mathcal{P}_r^k,n\mathcal{P}_s^k)=R(m\mathcal{P}_r^k,n\mathcal{C}_s^k)=((k-1)r+1)m+\Big\lfloor\frac{s+1}{2}\Big\rfloor n-1,$$
$$R(m\hat{\mathcal{P}}_r^k,n\hat{\mathcal{P}}_s^k)=(r+k-1)m+
(1+\lfloor\frac{s-1}{k}\rfloor)n-1,$$
$$R(m\mathcal{C}_r^k,n\mathcal{C}_s^k)=(k-1)rm+\Big\lfloor\frac{s+1}{2}\Big\rfloor n-1,$$
%$$R(m\hat{\mathcal{C}}_{kr}^k,n\hat{\mathcal{C}}_{ks}^k)=krm+sn-1,$$
and if $r>s$
$$R(m\mathcal{C}_r^k,n\mathcal{P}_s^k)=(k-1)rm+\Big\lfloor\frac{s+1}{2}\Big\rfloor n-1.$$

\end{conjecture}

Using Theorem \ref{goodramsey}, for a natural number $m$, a $k$-uniform hypergraph $\mathcal{H}$, a given $\mathcal{G}\in \cal{G}$$_k$ and sufficiently large $n$ we have

$$R(m\mathcal{G},n\mathcal{H})=m|V(\mathcal{G})|+n|V(\mathcal{H})|-m\alpha(\mathcal{G})-
1.$$

Based on this equality, we pose the following
conjecture for the Ramsey number of $k$-uniform hypergraphs corresponding to trees.

\medskip
\begin{conjecture} \label{trees}
For $k\geq 3$, assume that $\mathcal{G}=\mathcal{H}_k(T)$ and $\mathcal{H}=\mathcal{H}_k(T')$ where $T$ and $T'$ are trees in $\cal{G}$$_2$. If $\alpha(\mathcal{G})\leq \alpha(\mathcal{H})$ and $m\leq n$, then
$$R(m\mathcal{G},n\mathcal{H})=m|V(\mathcal{G})|+n|V(\mathcal{H})|-m\alpha(\mathcal{G})-
1.$$
\end{conjecture}

Using Theorems \ref{ramsey number for k-psths} and \ref{ramsey number for multiple 3-paths},  Conjecture \ref{trees} is true for an arbitrary $k$ when $T$ and $T'$ are either two paths with at most 4 edges or stars with $3$ edges and also for $k=3$ when
 $T$ and $T'$ are two arbitrary paths.
Clearly by Theorem \ref{main theorem}, if this conjecture holds for the case $m=n=1$, then it holds for every $m\geq n\geq1$.

\end{document}